\documentclass[12pt]{article}
\usepackage{amsmath,amsfonts,amssymb}

\newcommand{\C}{\mathbb C}
\newcommand{\Q}{\mathbb Q}
\newcommand{\Z}{\mathbb Z}
\newcommand{\cB}{\mathcal B}
\newcommand{\E}{\mathcal E}

\newcommand{\PP}{\mathbb P}
\newcommand{\Oh}{\mathcal O}
\newcommand{\eps}{\varepsilon}

\DeclareMathOperator{\rank}{rank}

\DeclareMathOperator{\Pic}{Pic}

\DeclareMathOperator{\hcf}{hcf}
\DeclareMathOperator{\lcm}{lcm}

\newtheorem{ingredient}{Ingredient}
\newtheorem{thm}{Theorem}
\newtheorem{definition}[thm]{Definition}
\newtheorem{notn}[thm]{Notation}

\newcommand{\QED}{\ifhmode\unskip\nobreak\fi\quad\ensuremath{\mathrm{QED}}}
\newcommand\st{{\ \vert\ }}

\title{Fano 3-folds of index 2}
\author{Gavin Brown \and Kaori Suzuki}
\date{}

\begin{document}
\maketitle

\begin{abstract}
We study Fano $3$-folds with Fano index $2$: that is, $3$-folds
$X$ with $\rank\Pic(X)=1$, $\Q$-factorial terminal singularities and
$-K_X=2A$ for an ample Weil divisor $A$.  We give a first
classification of all possible Hilbert series of such polarised
varieties $X,A$ and deduce both the nonvanishing of $H^0(X,-K_X)$
and the sharp bound $(-K_X)^3\ge 8/165$.
We list families that can be realised in codimension up to 4.
\end{abstract}

\section{Introduction}

We work over the complex number field $\C$ throughout,
and we denote the Picard number of $X$ by $\rho(X) = \rank\Pic(X)$.

\begin{definition}
A normal projective 3-fold $X$ is called a {\em Fano 3-fold} if and
only if $X$ has $\Q$-factorial terminal singularities, $-K_X$ is ample,
and $\rho(X) = 1$.

The {\em Fano index $f=f(X)$} of a Fano 3-fold $X$ is
\[
f(X) = \max \{ m \in \Z_{> 0} \st -K_X = mA
\textrm{ for some Weil divisor }A \}.
\]
where equality of divisors denotes linear equivalence of some multiple.
A Weil divisor $A$ for which $-K_X = fA$ is called
a {\em primitive ample divisor}.
\end{definition} 
Fano 3-folds are sometimes called $\Q$-Fano $3$-folds to distinguish
them from the classical nonsingular case.
By \cite{Su}, Theorem 0.3, we know that $f\le 19$.
The case $f=1$ is the main case with several hundred confirmed families
and much work towards classification ongoing.
In the case $f\ge 4$, \cite{Su1} contains a near classification of
about 80 families (only a handful of which remain in doubt).
We study the case $f=2$ here.
Typical classical examples
include the cubic 3-fold in $\PP^4$ and the intersection of two
quadrics in $\PP^5$;  see \cite{encyc}, Table~12.2, for example,
in which $r$ denotes the Fano index.
Note that the anti-canonical divisor $-K_X$ is only expected to be
a Weil divisor, although necessarily it will be $\Q$-Cartier.

A Fano $3$-fold $X$ with primitive ample divisor $A$ has a graded ring
\[
R(X, A) = \bigoplus_{n \geq 0} H^{0}(X, \Oh_{X}(nA)).
\] 
This graded ring is finitely generated.  The {\em Hilbert series of $X,A$}
is defined to be that of the graded ring $R(X,A)$.  A choice of minimal
(homogeneous) generating set $x_0,\dots,x_N\in R(X,A)$ determines
an embedding
\begin{equation}
\label{eq!emb}
X\hookrightarrow \PP^N = \PP(a_0,\dots,a_N)
\end{equation}
for some weighted projective space (wps) $\PP^N$, where
$x_i\in H^{0}(X, \Oh_{X}(a_iA))$.
With this embedding in mind, we say that $X,A$ has {\em codimension} $N-3$.

A Fano 3-fold is {\em Bogomolov--Kawamata stable},
or simply {\em stable}, if
\begin{equation*}
(-K_X)^3 \leq 3(-K_X c_2(X)).
\end{equation*}
Fano 3-folds satisfying this condition fall into one half of Kawamata's
argument on the boundedness of Fanos \cite{Ka2}.

Our first result is a list of the possible Hilbert series of graded rings
$R(X,A)$ for Fano $3$-fold of Fano index $2$ and their analysis in low
codimension.
\begin{thm}
\label{thm!main}
There are at most 1492 power series that are the Hilbert series
of some $X,A$ with $X$ a Fano 3-fold of Fano index~2 and $A$ a primitive
ample divisor, of which 1413 could correspond to stable Fano 3-folds.
Of these power series, 36 can be realised by some stable $X,A$ of
codimension $\le 3$.
\end{thm}

The proof is explained in section~\ref{sec!pf}; arguments of Kawamata
\cite{Ka1}, \cite{Ka2} and Suzuki \cite{Su} impose conditions
on geometrical data of $X,A$ which are then analysed by computer.
The 36 cases in low codimension are listed in
Tables~\ref{tab!cod1}--\ref{tab!cod3} in section~\ref{sec!cod123}.

In the case $f=1$, there are many Fano 3-folds with with empty
anti-canonical linear system.  This does not happen in index~2:
\begin{thm}
\label{thm!nonvan}
If $X,A$ is a Fano 3-fold of Fano index $2$
then $H^0(X,-K_X)\not= 0$.
\end{thm}

We showed in \cite{BS3} that $-K_X$ also has sections whenever
$f\ge 3$, and so this result shows that $-K_X$ has a section
whenever $f\not=1$.

When $f\ge 3$, the Riemann--Roch formula is a function of
the singularities of $X$.  Fano 3-folds with $f\le 2$ have
an extra parameter, their {\em genus}, $g=\dim H^0(X,A)-2$.
We have sharp bounds on the degree for various genera:
\begin{thm}
\label{thm!bounds}
Let $X,A$ be a Fano 3-fold of Fano index $2$ with $-K_X = 2A$.
Let $g=\dim H^0(X,A)-2$.
Then $-2\le g \le 9$, and the degree $A^3$ of $X$ is bounded below
by $1/165$; according to small $g$ the sharp bounds are:
\[
\begin{array}{r|ccccc}
\text{genus $g$} & -2 & -1 & 0 & 1 & 2  \\
\hline
\text{lower bound for $A^3$} & 1/165 & 1/35 & 1/3 & 1 & 2
\end{array}
\]
Moreover, all these lower bounds are achieved by stable Fano 3-folds.
\end{thm}

The graded ring approach to building classifications is well known;
we describe it in this case in section~\ref{sec!pf}.
There are two main points to be aware of.
First, although the the Hilbert series of any $R(X,A)$ will indeed
appear in the list of Theorem~\ref{thm!main}, there may be power
series in the list that do not correspond to such a graded ring:
the list comprises {\em candidates} for the Hilbert series of
Fano 3-folds, and an appearance on the list does not imply that a
Fano 3-fold exists with that Hilbert series.  
Second, we cannot say exactly which rings occur---for there will be
many degenerations, as in \cite{B2}, and classifying them will be
difficult.  But in many
cases we can predict at least the smallest possible codimension
of such a ring together with the weights of the corresponding wps.
When the proposed ring has codimension $\le 3$, then we are able
to construct it.  Gorenstein rings in codimension~4 are more subtle.
In section~\ref{sec!cod4} we describe a first classification into 35 cases.

In section~\ref{sec!examples} we discuss some typical examples.
Then section~\ref{sec!series} contains the proofs after first
assembling the ingredients:  the appropriate Riemann--Roch formula,
bounds on the singularities and degree, and Kawamata's boundedness result.
In section~\ref{sec!lists} we study the converse question, that of
constructing Fano 3-folds with given Hilbert series, listing results in
Tables~\ref{tab!cod1}--\ref{tab!cod4}.

It is our pleasure to thank Miles Reid for his help and
encouragement throughout this work.

\section{Examples}
\label{sec!examples}

We work out some examples, pointing out where each Ingredient~1--5
of section~\ref{sec!series} fits in.
We use our standard notation of Notation~\ref{stdnotn} below.
To run the graded ring method, we need a basket of singularities
and an integer $g\ge -2$.  In the first example, we choose $g=-2$.

Consider the basket of quotient singularities
\[
\cB =
\left\{
\frac{1}{3}(1,2,2),
\frac{1}{5}(1,4,2),
\frac{1}{11}(3,8,2)
\right\},
\]
where $\frac{1}{r}(a,-a,2)$ denotes the germ of the quotient
singularity $\C^3/(\Z/r)$ by the action
\[
(x,y,z) \mapsto (\eps^a x,\eps^{-a} y, \eps^2 z)
\quad
\textrm{where $\eps=\exp(2\pi i/r)$}.
\]
(The third component 2 is forced since such a singularity
will be polarised locally by the canonical class.)

We want to make a Fano 3-fold $X,A$ of Fano index~2 such
that $X$ has exactly the singularities of $\cB$
and $\dim H^0(X,A) = g + 2$.
(It would be enough that $X$ has terminal singularities
which contribute to Riemann--Roch as though they were the
singularities of $\cB$, but in practice a general element of
any family we construct has quotient singularities.)
According to Ingredient~\ref{ing!basket} in section~\ref{sec!series}
below, the basket $\cB$ is a possibility since
\[
\sum_{\cB} r - 1/r =
(3-1/3) + (5-1/5) + (11-1/11) < 24.
\]
Applying the formula of Ingredient~\ref{ingred!hs} with data $\cB,g$
describes Hilbert series
\[
P(t) = 
\frac{1-t^{38}}{(1-t^2)(1-t^3)(1-t^5)(1-t^{11})(1-t^{19})}.
\]
Certainly any hypersurface of the form
\[
X = X_{38} \subset \PP(2,3,5,11,19)
\]
will have Hilbert series $P_X(t)$ equal to $P(t)$.
One can check that a general member of this family is
indeed a quasismooth index~2 Fano 3-fold with singularities exactly
$\cB$ and $H^0(X,A)=0$.
This is one of eight hypersurfaces that are Fano 3-folds of
Fano index~2:  all eight are listed in Table~\ref{tab!cod1}
of section~\ref{sec!cod123}.

Notice that any quasismooth Fano 3-fold of this form has a K3 elephant;
that is, the section $S\in \mid{-}K_X\mid$ is a K3 surface,
$S_{38}\subset \PP(3,5,11,19)$.

\paragraph{Examples from del Pezzo surfaces}

A classical example is the Fano 3-fold $X\subset\PP^6$
whose hyperplane section is the del Pezzo surface of degree~5.
The equations of a general such 3-fold are well known:  they are the
five maximal Pfaffians of a skew $5\times 5$ matrix of general
linear forms on $\PP^6$.
Extensions of other del Pezzo surfaces also give rise
to index~2 Fano 3-folds.  In \cite{cascade}, such series of
families is called a `cascade' and other examples using log
del Pezzo surfaces are described there.

The main case of \cite{cascade} is a series of families which,
in order of increasing codimension (although, as with blowups
of $\PP^2$, the opposite order is also natural) begins
with the hypersurface $X_{10} \subset \PP(1^2, 2, 3, 5)$,
followed by $X_{4,4} \subset \PP(1^3, 2^2, 3)$ in codimension~2
and a family in codimension~3, $X\subset \PP(1^4, 2^2, 3)$.

These families arise using $\cB=\{\frac{1}{3}(1,2,2)\}$
and allowing $g$ to vary, starting at $g=0$ and increasing.
When $g=0$, the formula of Ingredient~\ref{ingred!hs} with
data $\cB,g$ describes Hilbert series
\[
P(t) = (1-t^{10})/(1-t)^2(1-t^2)(1-t^3)(1-t^5),
\]
and, as before, this corresponds to a family of Fano hypersurfaces
above of degree $A^3=1/3$.  When $g=1$, we get
\[
P(t) = (1-t^4)^2/(1-t)^3(1-t^2)^3(1-t^3),
\]
corresponding to the codimension~2 complete intersection above
of degree $4/3=1/3+1$.
Increasing $g$ again to $g=2$, and thus $A^3=7/3=2/3+2$ by
Ingredient~\ref{ing!mindeg} below, describes Hilbert series
\[
P(t) = (t^4 + t^3 + 3t^2 + t + 1)/(1-t)^3(1-t^3).
\]
This description is not as revealing as previous ones, but it is
easy to remedy.  The numerator indicates that a further generator
in degree~1 is needed, and with an extra factor of $1-t$, we see that
\[
P(t) = (-t^5 - 2t^3 + 2t^2 + 1)/(1-t)^4(1-t^3).
\]
Now it is clear that two new generators in degree~2 are
required---these also serve to polarise the index~3 singularity in
the basket.  The final result is
\[
P(t) = (-t^9 + 2t^6 + 3t^5 - 3t^4 - 2t^3 + 1)/(1-t)^4(1-t^2)^2(1-t^3).
\]
This Hilbert series is realised by a graded ring with generators
$x_{1}, x_{2}, x_{3}$, $y_{1}, y_{2}$, $z$ in degrees 1, 1, 1, 2, 2, 3 
and relations generated by the five maximal Pfaffians
of the following $5 \times 5$ skew matrix (where as usual we
omit the leading diagonal of zeroes and leave the skew lower half
implicit)
\[
M = \begin{pmatrix}
      x_{1} & x_{2} & b_{14} & b_{15}\\
      & x_{3} & b_{24} & b_{25}\\ 
      &&  b_{34} & b_{35}\\
      &&&   z
   \end{pmatrix}
\quad
\textrm{of degrees}
\quad
   \begin{pmatrix} 
     1 & 1 & 2 & 2\\
     & 1 & 2 & 2\\ 
     &&  2 & 2\\
     &&&   3
   \end{pmatrix}.
\]
The $b_{ij}$ are forms of degree~2 in the $x$ and $y$ variables.
It is easy to check that for general $b_{ij}$, this defines the
graded ring of a Fano 3-fold in $\PP^6(1^4,2^2,3)$ as required, and
we think of the $b_{ij}$ as being parameters defining a flat family
of Fano 3-folds of index~2 with given basket $\cB$ and genus $g$.

\section{Hilbert series of Fano 3-folds}
\label{sec!series}

We describe in sections~\ref{sec!RR}--\ref{sec!bddness} the five
ingredients that go into the raw list of power series that is the
basis for our classification.  Our raw list will include the Hilbert
series of every index~2 Fano 3-fold, although it may contain power
series not of this form.  The derivation of the list is explained
in section~\ref{sec!pf}.
Although it is not necessary for everything below, we suppose throughout
that Fano 3-folds $X$ treated here will have Fano index $f(X)=2$.

\subsection{Riemann-Roch formula}
\label{sec!RR}

We explain the notion of a {\em basket of singularities}; see
\cite{YPG} or \cite{B} section~2.1.
Let $P\in X$ be a 3-dimensional terminal singularity of index $r > 1$
and $(U, P)$ a germ at $P$.  (In particular, $rK_U$ is Cartier on $U$.)

If $(U, P)$ is a quotient singularity, then it isomorphic to the germ
at the origin of some $\C^3/(\Z/r)$, with action
\[
\eps\cdot(x,y,z)\mapsto (\eps^ax, \eps^ay, \eps^cy).
\]
We denote this by $\frac{1}{r}(a, b, c)$
and recall that $\hcf(r,abc) = 1$ and $a+b=r$ up to permutations of $a,b,c$.
When $X$ has Fano index~2, then we may suppose that
$(a,b,c) = (a,r-a,2)$.  In particular, $r$ cannot be even.


If $(U, P)$ is not a quotient singularity, then
$(U, P)$ can be deformed to a unique finite collection of
terminal quotient singularities, say $\{Q_1,\dots, Q_{n(P)}\}$ (where,
as is usual, this is a set with possible repetitions) for some
number $n(P) \geq 1$.  Each point $Q_i$ is some quotient singularity
$\frac{1}{r_{i}}(a_i, -a_i,2)$
where $r_i$ and $a_i$ are coprime,
$r_i \geq 2$, and $r=\lcm\{r_1,\dots r_{n(P)}\}$.
We call the set $\cB(U, P) := \{Q_{1},\dots Q_{n(P)}\}$ the
{\em basket of singularities} of $(U, P)$. 

In the global case, we assemble all local baskets into one.
\begin{definition}
Let $X$ be a 3-fold with terminal singularities and $\{P_1,\dots,P_m\}$
the set of singular points of $X$ of index $\geq 2$.
Denoting germs $P_{i} \in U_{i} (\subset X)$,
we define the {\em basket of singularities of $X$}
to be the disjoint union $\cB(U_{1},P_{1}) \cup\dots\cup \cB(U_m,P_m)$
(a set with possible repetitions.)
\end{definition} 

\begin{notn}
\label{stdnotn}
We denote a Fano 3-fold of index~$f=2$ for which $-K_X=2A$ by $X,A$.
The basket of singularities of $X,A$ is denoted $\cB$.
A typical singularity of $\cB$ is denoted $\frac{1}{r}(a,-a,2)$,
and we use this notation whenever taking a sum over the elements of $\cB$.
\end{notn}

\begin{thm}[\cite{Ka1}]
Let $X$ be a Fano $3$-fold with basket of singularities $\cB$.
Then $-K_Xc_2(X) > 0$, the Euler characteristic $\chi(\mathcal{O}_X)=1$
and
\[
\chi(\mathcal{O}_{X}) =
\frac{-K_{X}c_{2}(X)}{24} + \sum\frac{r^{2} -1}{24r},
\] 
the sum taken over $\cB$ (see Notation~\ref{stdnotn}).
\end{thm}
This theorem gives the following bounds on the
singularities in the basket:
\begin{ingredient}[Basket bound]
\label{ing!basket}
If $\cB$ is the basket of a Fano 3-fold $X$, then
\[
\sum\frac{r^{2} -1}{r} < 24
\quad
\textrm{and}
\quad
-K_{X}c_{2}(X) = 24 - \sum\frac{r^{2} -1}{r}
\] 
where each sum is taken over $\cB$ (see Notation~\ref{stdnotn}).
\end{ingredient}

For the next theorem, recall that the plurigenera of a polarised
variety $X,A$ are denoted $P_n(X,A) = \dim H^0(X,nA)$.

\begin{thm}[{\cite{Su}}]
\label{thm!RR}
Let $X$ be a Fano $3$-fold of Fano index $f$ and basket of singularities
$\cB$.  Let $A$ be a primitive Weil divisor with $-K_{X} = fA$.  Then 
\begin{eqnarray*}
\chi(\mathcal{O}_{X}(nA)) & = &
\chi(\mathcal{O}_{X}) + \frac{n(n+f)(2n+f)}{12}A^{3}
 + \frac{nAc_{2}(X)}{12} \\
 &&
 + \sum_\cB \left( -i_n\frac{r^{2} -1}{12r} 
 + \sum_{j=1}^{i_n-1} \frac{\overline{bj}(r-\overline{bj}) }{2r}\right)
\end{eqnarray*}
where the sum is taken over points $P=\frac{1}{r}(a,-a,2)$
in $\cB$ (see Notation~\ref{stdnotn}) using notation:
$i_n \in [0, r-1]$ is the local index of $nA$ at $P$ (see \cite{YPG}),
$b \in [0, r-1]$ satisfies $ab \equiv 2 \mod r$ and 
$\overline{c} \in [0, r-1]$ is the residue of $c$ $\mod r$.

By Kawamata--Viehweg vanishing, $\chi(nA) = h^{0}(nA)$ for all $n > -f$.
So the Hilbert series $P_{X,A}(t)$ of $X$ is:
\begin{eqnarray*}
P_{X,A}(t) & = & \sum_{n = 0}^{\infty} P_{n}(X,A)t^{n}\hspace*{1cm}\\
 & = &
  \frac{1}{1-t} +
  \frac{(f^{2}+3f+2)t+(-2f^{2}+ 8)t^{2}+(f^{2}-3f + 2)t^{3}}{12(1-t)^{4}}A^{3}\\
 & & \qquad + \frac{t}{(1-t)^{2}}\frac{Ac_{2}(X)}{12} + \sum_{P\in\cB} c_P(t)
\end{eqnarray*}
where, for a point $P=\frac{1}{r}(a,-a,2)$ in $\cB$,
\begin{eqnarray*}
c_P(t) & = & \frac{1}{1-t^{r}}
  \left(
  \sum_{k=1}^{r-1}
    \left(-i_k\frac{r^2-1}{12r} +
    \sum_{j=1}^{i_k-1}\frac{\overline{bj}(r-\overline{bj})}{2r}\right)t^k
  \right).
\end{eqnarray*}
\end{thm}
Setting $f=2$ in the Hilbert series above gives the closed formula:
\begin{ingredient}[Hilbert series]
Using Notation~\ref{stdnotn} and the expression for $c_P(t)$
from Theorem~\ref{thm!RR}, the Hilbert series is
\label{ingred!hs}
\[
P_{X,A}(t) = \frac{1}{1-t} + \frac{t}{(1-t)^{4}}A^{3}
  + \frac{t}{(1-t)^{2}}\frac{Ac_{2}(X)}{12} + \sum c_P(t)
\]
where the sum is taken over $\cB$.
\end{ingredient}
Setting $n=1$ in the Riemann--Roch formula, we compute the minimum
possible value of the degree $A^3$.
\begin{ingredient}[Minimum degree]
\label{ing!mindeg}
Using Notation~\ref{stdnotn},
\[
0 < A^3 =
 - 1 - 
  \frac{Ac_2(X)}{12} -
    \sum \left( -i_1\frac{r^{2} -1}{12r} 
     + \sum_{j=1}^{i_1-1} \frac{\overline{bj}(r-\overline{bj}) }{2r}\right)
 + N
\]
for some integer $N\ge 0$, and where the sum is taken over $\cB$.

\end{ingredient}
Since the Riemann--Roch formula also holds when $n=-1$, we have an
additional constraint.
\begin{ingredient}[Polarisation condition]
\label{ing!pol}
Using Notation~\ref{stdnotn},
\[
1 +
\sum \left( -(i_{-1})\frac{r^{2} -1}{12r} 
 + \sum_{j=1}^{(i_{-1})-1} \frac{\overline{bj}(r-\overline{bj}) }{2r}\right) =
\frac{Ac_{2}(X)}{12}.
\]
where the sum is taken over $\cB$.
\end{ingredient}

\subsection{Kawamata Boundedness}
\label{sec!bddness}

Following \cite{Ka2}, let $X$ be a Fano 3-fold with
Fano index $2$ and  $\E = (\Omega_{X}^{1})^{**}$ the double dual of
the sheaf of K\"ahler differentials of $X$. 
We do not define $\mu$-semistability (with respect to $-K_{X}$) here,
since we do not use it further, but note the role it plays in
strengthening the following bound on the degree.
\begin{thm}[{\cite{Ka2}}]
\label{thm:kaw}
Let $X$ be a Fano $3$-fold with Fano index $2$.
Then $(-K_{X})^{3} \leq \frac{16}{5}(-K_{X}c_{2}(X))$. 
If, furthermore, the sheaf $\E$ described above is $\mu$-semistable,
then $(-K_{X})^{3} \leq 3(-K_{X}c_{2}(X))$.
\end{thm} 

\begin{ingredient}
\label{ing!kawa}
Using Notation~\ref{stdnotn}, $A^3 \le \frac{4}{5}Ac_2(X)$.
(In the stable case, the upper bound is $A^3\le \frac{3}{4}Ac_2(X)$.)
\end{ingredient}

Applying the upper bound for $Ac_2$ of Ingredient~\ref{ing!basket},
this bound implies $A^3\le 48/5$ (and $A^3\le 9$ in the stable case.)
In fact, once all possible baskets are calculated, the sharp bound
is $A^3\le 9$ achieved using the empty basket and $g=9$;  the next
largest degree is $A^3=25/3$ achieved by the basket $\{\frac{1}{3}(1,2,2)\}$
and $g=8$ (which is not stable since $Ac_2=32/3$).
In particular, for a fixed basket there are at
most 9 different values for $g$ that give Fano Hilbert series.
Comparing with \cite{cascade}, we regard this as a bound on the
number of blowups (or projections) that we can make from maximal $g$.
A familiar instance of this bound is the maximum number of blowups of
$\PP^2$ that is a del Pezzo surface;  again, 8 is the limit.
In \cite{cascade}, such a cascade of 8 log del Pezzo surfaces is
constructed linking the hypersurface $S_{10}\subset\PP(1,2,3,5)$ with
a log del Pezzo of degree $25/3$.

\subsection{Proofs of Theorems~\ref{thm!main}--\ref{thm!bounds}}
\label{sec!pf}
These proofs use some computer calculations which we describe rather
than reproducing; we use the computer algebra system Magma \cite{M}
for our calculations, and a short file with Magma code that
can be run to generate these results is available at \cite{code}.

\paragraph{Listing the Hilbert series}

The first step is to construct the list of all power series according
to the five ingredients assembled in section~\ref{sec!series}.
This requires little comment.  We simply compute all possible baskets
satisfying Ingredients~\ref{ing!basket} and \ref{ing!pol}, together
with all possible values for $N$ using Ingredients~\ref{ing!mindeg}
and \ref{ing!kawa}, and then apply the Riemann--Roch formula
of Ingredient~\ref{ingred!hs}.
The result is all power series that could be the Hilbert series
of a Fano 3-fold of index~2.  There are 1492 such series in all.
(If we impose the lower `stable' bound of Ingredient~5, this
number reduces to 1413.)

\paragraph{Bounds on the degree}
It is easy to compute both upper and lower bounds on the degree
(and also on $Ac_2(X)$) by computer check on the list.
In Theorem~\ref{thm!bounds} we list lower bounds for small $g$, since
our methods of construction reveal Fano 3-folds that realise these
bounds.  The following table includes all the degree bounds for each
genus $g$, even though we do not know whether all of them are sharp
or not.
\begin{center}
$
\begin{array}{|r|ccccc}
\hline
\text{genus $g$} & -2 & -1 & 0 & 1 & 2  \\
\hline
\text{lower bound for $A^3$} & 1/165 & 1/35 & 1/3 & 1 & 2\\
\text{upper bound for $A^3$} & 11/15 & 32/21 & 89/39 & 64/21 & 19/5\\
\hline
\text{total number of series} & 337& 470& 303& 174& 97\\
\text{(of which are unstable)} & 6& 14& 14& 15& 11\\
\hline
\end{array}
$
\hfill
\vspace{3mm}
\\ \hfill
$
\begin{array}{ccccccc|}
\hline
3 & 4 & 5 & 6 & 7 & 8 & 9 \\
\hline
3 & 4 & 5 & 6 & 7 & 8 & 9 \\
68/15 & 16/3 & 6 & 48/7 & 38/5 & 25/3 & 9\\
\hline
54& 28& 14& 8& 4& 2& 1\\
7& 5& 2& 2& 2& 1& 0\\
\hline
\end{array}
$
\end{center}

\paragraph{Nonempty anticanonical system}
The proof of Theorem~\ref{thm!nonvan} is straightforward: since we
do not discard candidate Hilbert series unless they are proved not
to come from a Fano 3-fold, it is enough to calculate
the coefficient of $t^2$ in each one.  Again, this is done by
a computer on all 1492 Hilbert series.

Since $H^0(X,-K_X)\not= 0$, we can ask whether the linear system
$\vert -K_X \vert$ contains a K3 surface.  Such a K3 section is sometimes
impossible because its singular rank would be too big---see \cite{B},
Proposition~4.  For example, $\cB=\{\frac{1}{21}(10,11,2)\}$ with $g=0$
satisfy all our numerical conditions (including stability)---they predict
$A^3 = 19/21$ and $(1/12)Ac_2(X) = 8/63$.  But a K3 section cannot exist
because the corresponding surface singularity $\frac{1}{21}(10,11)$ has
20 exceptional curves in its resolution, pairwise orthogonal in
the Picard group, which cannot happen in $H^2(S)$ for a K3 surface $S$.
There are 171 such cases (of which 9 are unstable).  (These examples
all appear to be in high codimension, and we do not analyse them further.)


When the singular rank is $\le 19$, we get an estimate of the degrees
of generators of a model $R(X,A)$ for a Fano 3-fold by comparing with
K3 surfaces appearing in the K3 database \cite{B}.  This gives an idea of
how to understand the result of Theorem~\ref{thm!main}, just as the estimates
of \cite{cascade} 3.2.5 do in the case $f\ge 3$.  But this list should not
be taken as more than a guide.  Using this K3 comparison, the number
of Hilbert series per codimension is as follows.
\begin{center}
$
\begin{array}{|r|cccccccc}
\hline
\textrm{estimated codimension} & 1&2&3&4&5&6&7&8\\
\hline
\text{total number of series} &
 8 & 26 & 2 & 35 & 13 & 59 & 25 & 99\\
\text{(of which are unstable)} &
 0 & 0 & 0 & 0 & 0 & 0 & 0 & 0\\
\hline
\end{array}
$
\hfill
\vspace{3mm}
\\ \hfill
$
\begin{array}{cccccccccc|}
\hline
9&10&11&12&13&14&15&16&17&18\\
\hline
51 & 163 & 93 & 227 & 126 & 255 & 48 & 78 & 8 & 3 \\
1 & 2 & 4 & 6 & 8 & 30 & 6 & 11 & 0 & 0\\
\hline
\end{array}
$
\end{center}
We describe the 36 candidates in codimension~$\le 3$ in
section~\ref{sec!cod123} and the 35 candidates in codimension~4
in section~\ref{sec!cod4} below.

\section{Constructing lists of Fano $3$-folds}
\label{sec!lists}

The remaining claim is that 36 of the 1492 families can be
realised in codimension $\le 3$.  This is achieved by the
examples in Tables~\ref{tab!cod1}--\ref{tab!cod3} below.
We also explain the list in Table~\ref{tab!cod4} of codimension~4
candidates and describe the role of unprojection methods.

\subsection{Models in low codimension}
\label{sec!results}

\paragraph{Finding the first generators}
The following standard analysis of a Hilbert series bounds the
minimum number of generators below.  Since we are concentrating
on codimension $\le 3$, once we have confirmed that the ring must
have at least eight generators if it exists, then we do not pursue
it further.

Consider the data $\cB = \{ \frac{1}{9}(1,8,2) \}$, $g=1$.
Suppose there is a Fano $X,A$ with this data;  we begin
to describe $R(X,A)$.
(Such $X$ would have $A^3=13/9$ and $Ac_2=17/27$.)
The Hilbert series expressed as a power series is
\[
1 + 3t + 8t^2 + 17t^3 + 32t^4 + 54t^5 + 85t^6 + \cdots.
\]
Any graded ring having this Hilbert series must have exactly
three generators in degree~1.  These span at most a 6-dimensional
subspace in degree~2, so there are at least two new generators
in degree~2.  If this number is exactly two, then a similar argument
shows that there must be at least one new generator in degree~3.
The only alternative is to have three new generators in degree~2.
Either way, we already have at least six generators in the ring.

Now we turn to the singularities.  Certainly $X$ has a terminal
singularity of index~9.  Whether it is a quotient singularity
or not, locally it must be a quotient by $\Z/9$ acting with
at least one eigencoordinate of weight 8---every case of index~9 in
the classification of terminal singularities \cite{mori} and
their baskets \cite{YPG} has local eigencoordinates $a,-a$,
where $a=1$ in this case.  So $R(X,A)$ must
have at least one generator of weight divisible by 9 and another
of weight 8 modulo 9.  So $R(X,A)$ is in codimension $\ge~4$.
(One could also argue on the denominator of the Hilbert series:
as a rational function, it has cyclotomic polynomial of degree~9
in its denominator, so the ring needs a generator in degree some
multiple of 9 to cancel this.)

\paragraph{Constructing varieties and additional generators}

We carry out the analysis above systematically using a computer.
When it suggests that a ring may be in codimension $\le 3$, then we
attempt to construct it as in the examples of section~\ref{sec!examples}.
The construction is straightforward, although there is one small twist:
it happens frequently that we can find a graded ring with the right
Hilbert series but that does not correspond to a Fano 3-fold.

For example, with $\cB = \{ \frac1{11}(2,9,2) \}$, $g=-1$,
the Hilbert series is
\[
P(t)=\frac{(1-t^6)(1-t^8)(1-t^{10})}{\prod(1-t^a)}
\]
where the product is taken over $a\in\{1,2,2,2,3,5,11\}$.
This is the Hilbert series of any $X_{6,9,10}\subset\PP(1,2,2,3,3,5,11)$.
But such a variety is not a Fano 3-fold: the equations cannot involve
the variable of degree 11, so it will have a cone singularity at
the index~11 point.
Moreover, there is no variable of degree 9 to polarise the singularity
as we wanted.  In this case, the solution is clear:
adding a generator of degree~9 to the ring gives a codimension~4 model
$X\subset\PP(1,2,2,2,3,5,9,11)$.
The family of complete intersections comprise a component of the
Hilbert scheme that does not contain the desired varieties.


\subsection{Classification in codimensions 1 to 3}
\label{sec!cod123}

The lists in Tables~\ref{tab!cod1}--\ref{tab!cod3} contain families
whose general element lies in
codimension at most 3.  (We denote the singularity $\frac{1}{3}(1,2,2)$
simply by $\frac{1}{3}$.)  In addition to these, there will be
degenerations of low codimension families that occur in higher
codimension.  For example, a Fano $X_{6,18}\subset\PP(1,2,3,5,6,9)$
for which the variable of degree~6 does not appear in the equation
of degree~6 is not listed in the table of codimension~2 Fano
3-folds;  its Hilbert series is that of the degree~$1/15$ hypersurface,
and the degeneration occurs as the two singularities of index~3
in that hypersurface come together.

It is easy to confirm that the general element in each case 
of Tables~\ref{tab!cod1}--\ref{tab!cod3} is a Fano 3-fold with
the indicated properties.

\begin{table}[h]
\[
\begin{array}{lccc}
\quad\textrm{Fano hypersurface} & \textrm{Basket $\cB$} & A^3 & Ac_2/12 \\
\hline
\hline
X_3\subset\PP(1,1,1,1,1) & \textrm{nonsingular} & 3 & 1 \\
X_4\subset\PP(1,1,1,1,2) & \textrm{nonsingular} & 2 & 1 \\
X_6\subset\PP(1,1,1,2,3) & \textrm{nonsingular} & 1 & 1 \\
X_{10}\subset\PP(1,1,2,3,5) & \frac{1}{3} \textrm{ ($= \frac{1}{3}(1,2,2)$)} & 1/3 & 8/9 \\
\hline
X_{18}\subset\PP(1,2,3,5,9) &
  2\times\frac{1}{3},\ \frac{1}{5}(1,4,2) & 1/15 & 26/45 \\
X_{22}\subset\PP(1,2,3,7,11) &
  \frac{1}{3},\ \frac{1}{7}(3,4,2) & 1/21 & 38/63 \\
X_{26}\subset\PP(1,2,5,7,13) &
  \frac{1}{5}(2,3,2),\ \frac{1}{7}(1,6,2) & 1/35 & 18/35 \\
X_{38}\subset\PP(2,3,5,11,19) &
  \quad\frac{1}{3},\ \frac{1}{5}(1,4,2),\ \frac{1}{11}(3,8,2)\quad & 1/165 & 116/495\\
\hline
\end{array}
\]
\caption{Fano 3-folds in codimension 1\label{tab!cod1}}
\end{table}

\begin{table}[t]
\[
\begin{array}{lccc}
\qquad\textrm{Fano 3-fold} & \textrm{Basket $\cB$} & A^3 & Ac_2/12 \\
\hline
\hline
X_{2,2}\subset\PP(1,1,1,1,1,1) &
  \textrm{nonsingular}
  & 3 & 1 \\
X_{4,4}\subset\PP(1,1,1,2,2,3) &
  \frac{1}{3}
  & 4/3 & 8/9 \\
X_{4,6}\subset\PP(1,1,2,2,3,3) &
  2\times \frac{1}{3}
  & 2/3 & 7/9 \\
X_{6,6}\subset\PP(1,1,2,2,3,5) &
  \frac{1}{5}(2,3,2)
  & 3/5 & 4/5 \\
X_{6,8}\subset\PP(1,1,2,3,4,5) &
  \frac{1}{5}(1,4,2)
  & 2/5 & 4/5 \\
\hline
X_{6,6}\subset\PP(1,2,2,3,3,3) &
  4\times\frac{1}{3}
  & 1/3 & 5/9 \\
X_{6,8}\subset\PP(1,2,2,3,3,5) &
  2\times\frac{1}{3},\ \frac{1}{5}(2,3,2)
  & 4/15 & 26/45 \\
X_{6,10}\subset\PP(1,2,2,3,5,5) &
  2\times \frac{1}{5}(2,3,2)
  & 1/5 & 3/5 \\
X_{8,10}\subset\PP(1,2,2,3,5,7) &
  \frac{1}{3},\ \frac{1}{7}(2,5,2)
  & 4/21 & 38/63 \\
X_{10,14}\subset\PP(1,2,2,5,7,9) &
  \frac{1}{9}(2,7,2)
  & 1/9 & 17/27 \\
\hline
X_{8,10}\subset\PP(1,2,3,4,5,5) &
  \frac{1}{3},\ 2\times\frac{1}{5}(1,4,2)
  & 2/15 & 22/45 \\
X_{8,12}\subset\PP(1,2,3,4,5,7) &
  \frac{1}{5}(1,4,2),\ \frac{1}{7}(3,4,2)
  & 4/35 & 18/35 \\
X_{10,12}\subset\PP(1,2,3,4,5,9) &
  \frac{1}{3},\ \frac{1}{9}(4,5,2)
  & 1/9 & 14/27 \\
X_{12,14}\subset\PP(1,2,3,4,7,11) &
  \frac{1}{11}(4,7,2)
  & 1/11 & 6/11 \\
X_{10,12}\subset\PP(1,2,3,5,6,7) &
  2\times\frac{1}{3},\ \frac{1}{7}(1,6,2)
  & 2/21 & 31/63 \\
X_{14,16}\subset\PP(1,2,5,7,8,9) &
  \frac{1}{5}(2,3,2),\ \frac{1}{9}(1,8,2)
  & 2/45 & 58/135 \\
\hline
X_{10,12}\subset\PP(2,2,3,5,5,7) &
  2\times\frac{1}{5}(2,3,2),\ \frac{1}{7}(2,5,2)
  & 2/35 & 11/35 \\
X_{10,14}\subset\PP(2,2,3,5,7,7) &
  \frac{1}{3},\ 2\times\frac{1}{7}(2,5,2)
  & 1/21 & 20/63 \\
X_{12,14}\subset\PP(2,2,3,5,7,9) &
  \frac{1}{3},\ \frac{1}{5}(2,3,2),\ \frac{1}{9}(2,7,2)
  & 2/45 & 43/135 \\
X_{14,18}\subset\PP(2,2,3,7,9,11) &
  2\times\frac{1}{3},\ \frac{1}{11}(2,9,2)
  & 1/33 & 32/99 \\
X_{18,22}\subset\PP(2,2,5,9,11,13) &
  \frac{1}{5}(1,4,2),\ \frac{1}{13}(2,11,2)
  & 1/65 & 17/65 \\
\hline
X_{10,12}\subset\PP(2,3,3,4,5,7) &
  4\times\frac{1}{3},\ \frac{1}{7}(3,4,2)
  & 1/21 & 17/63 \\
X_{12,14}\subset\PP(2,3,4,5,7,7) &
  \frac{1}{5}(2,3,2),\ 2\times\frac{1}{7}(3,4,2)
  & 1/35 & 8/35 \\
X_{14,16}\subset\PP(2,3,4,5,7,11) &
  \frac{1}{3},\ \frac{1}{5}(2,3,2),\ \frac{1}{11}(4,7,2)
  & 4/165 & 116/495 \\
X_{18,20}\subset\PP(2,4,5,7,9,13) &
  \frac{1}{7}(2,5,2),\ \frac{1}{13}(4,9,2)
  & 1/91 & 16/91 \\
X_{18,20}\subset\PP(2,5,6,7,9,11) &
  \frac{1}{3},\ \frac{1}{7}(2,5,2),\ \frac{1}{11}(5,6,2)
  & 2/231 & 103/693 \\
\hline
\end{array}
\]
\caption{Fano 3-folds in codimension 2}
\end{table}

\begin{table}
\[
\begin{array}{lccc}
\qquad\textrm{Fano 3-fold} & \textrm{Basket $\cB$} & A^3 & Ac_2/12 \\
\hline
\hline
X_{2,2,2,2,2}\subset\PP(1,1,1,1,1,1,1) & \textrm{\quad nonsingular\quad} & 5 & 1 \\
X_{3,3,4,4,4}\subset\PP(1,1,1,1,2,2,3) & \frac{1}{3} & 7/3 & 8/9 \\
\hline
\end{array}
\]
\caption{Fano 3-folds in codimension 3\label{tab!cod3}\label{tab!cod2}}
\end{table}

\subsection{Classification in codimension 4}
\label{sec!cod4}

We continue our analysis of rings into higher codimension,
although the method becomes more complicated and we do not
check all possibilities rigorously.  (Complete results are on the
webpage \cite{code}.)  We list in Table~\ref{tab!cod4} the 35
Hilbert series with proposals for models in codimension~4---it is
conceivable that there are other examples (other than degenerations)
but we do not expect them.  We list them according to the weights
of the ambient wps $\PP^7$, that is, the degrees of minimal
generators of $R(X,A)$; for full details, see the webpage at \cite{code}.
The task is to construct these Fano 3-folds.

\paragraph{Projection guided by a K3 section}

If $|{-}K_X|$ contains a K3 surface $S$ we may compare our results
with those of \cite{B} as a guide to the constructions we might make.
Consider, the following codimension~4 candidate:
\[
X\subset\PP(2,2,3,5,5,7,12,17)
\]
with $\cB=\{\frac{1}{17}(5,12,2)\}$.
A section in degree~2 is a known K3 surface known with a
Type~1 projection from a $\frac{1}{17}(5,12)$ point:
\[
\left(
S\subset\PP(2,3,5,5,7,12,17)
\right)
\dashrightarrow
\left(
T_{10,12,14,15,17}\subset\PP(2,3,5,5,7,12)
\right).
\]
The K3 surface $T$ can be constructed easily, and moreover it can
be forced to contain a linear $\PP(5,12)$.  This curve can be
unprojected to construct $S$.
This is how we proceed with cases in codimension~4 here, but there is
a twist.

When we project from a quotient singularity---say the point
$\frac{1}{17}(5,12,2)$ in the example above---the result will
always contain a line of index~2 singularities.  That is, projection
automatically incurs canonical singularities that contribute to
Riemann--Roch.  (Recall from \cite{ABR} that projection
of 3-folds results in slightly worse singularities that we
usually work with, but that typically these are Gorenstein and
do not contribute in the Riemann--Roch formula.)
And so the result of projection will not have Hilbert series
already on our lists.  This obstructs the inductive approach
exemplified by \cite{B}.

There are two ways out.  One would be to mimic classical
methods and make a double projection.  In the example above, after
projecting from $\frac{1}{17}(5,12,2)$, we can project from the resulting
$\frac{1}{12}(5,7,2)$ singularity.  This second projection contracts
the index~2 line, and we land in the family
\[
Z_{10,12}\subset\PP(2,2,3,5,5,7).
\]
To proceed, we would need an analysis of the exceptional locus of
such double projections and a Type~III-style unprojection
result, following \cite{kinosaki} e.g.~9.16.

The approach we take is to consider also weak Fano 3-folds with
canonical singularities.  This avoids the bottleneck in codimension~$3$
for general unprojection methods: by admitting some canonical
singularities, in particular lines of
index~2 with up to two non-isolated points of higher index on them,
we see many more weak Fano 3-folds in codimension~3.
We impose appropriate surfaces in these
3-folds and the apply general results of unprojection;
this method is explained in detail in \cite{BKR}.  In short, we
construct the example above by imposing the plane $\PP(2,5,12)$
as a linear subspace in the codimension~3 weak Fano 3-fold
\[
\PP(2,5,12)\subset Y_{10,12,14,15,17}\subset\PP(2,2,3,5,5,7,12).
\]
Such $Y$ has a non-isolated singularity $\frac{1}{12}(2,5,7)$ on
a line of index~2 singularities, which are all strictly canonical
singularities.

\paragraph{Codimension~4 candidates having no Type~I projection}
The projection method outlined above constructs examples for
33 of the 35 candidates in codimension~4.  The remaining two
cases, which we do not construct, are:
\[
X\subset\PP^7(2,3,3,4,5,5,6,7)
\]
with $\cB=\{ 5\times\frac{1}{3}(1,2,2), \frac{1}{5}(1,4,2)\}$,
$A^3 = 1/15$,  $\frac{1}{12}Ac_2(X) = 11/45$ and Hilbert numerator
$1 - t^8 - t^9 - 2t^{10} - t^{11} - 2t^{12} + t^{14} +\cdots+ t^{33}$;
and
\[
X\subset\PP^7(2,3,5,6,7,7,8,9)
\]
with
$\cB=\{ 3\times\frac{1}{3}(1,2,2), \frac{1}{5}(2,3,2), \frac{1}{7}(1,6,2)\}$,
$A^3 = 1/35$, $\frac{1}{12}Ac_2(X) = 19/105$ and Hilbert numerator
$1 - t^{12} - t^{13} - 2t^{14} - t^{15} - 2t^{16} - t^{17} - t^{18}
+ t^{19} +\cdots + t^{45}$.

\begin{table}
$
\begin{array}{lccc}
\quad\textrm{Ambient $\PP^7$} & \textrm{Basket $\cB$} & A^3 & Ac_2/12 \\
\hline
\hline
\PP( 1, 1, 1, 1, 1, 1, 1, 1 ) & \textrm{nonsingular}& 6& 1 \\
\PP( 1, 1, 1, 1, 1, 2, 2, 3 ) & \frac{1}{3}& 10/3& 8/9 \\
\PP( 1, 1, 1, 2, 2, 2, 3, 3 ) & 2\times\frac{1}{3}& 5/3& 7/9 \\
\PP( 1, 1, 1, 2, 2, 2, 3, 5 ) & \frac{1}{5}(2,3,2)& 8/5& 4/5 \\
\PP( 1, 1, 1, 2, 2, 3, 4, 5 ) & \frac{1}{5}(1,4,2)& 7/5& 4/5 \\
\hline
\PP( 1, 1, 2, 2, 2, 3, 3, 3 ) & 3\times\frac{1}{3}& 1& 2/3 \\
\PP( 1, 1, 2, 2, 2, 3, 3, 5 ) & \frac{1}{3}, \frac{1}{5}(2,3,2)& 14/15& 31/45 \\
\PP( 1, 1, 2, 2, 2, 3, 5, 7 ) & \frac{1}{7}(2,5,2)& 6/7& 5/7 \\
\PP( 1, 1, 2, 2, 3, 3, 4, 5 ) & \frac{1}{3},\frac{1}{5}(1,4,2)& 11/15& 31/45 \\
\PP( 1, 1, 2, 2, 3, 3, 4, 7 ) & \frac{1}{7}(3,4,2)& 5/7& 5/7 \\
\hline
\PP( 1, 1, 2, 3, 4, 5, 6, 7 ) & \frac{1}{7}(1,6,2)& 3/7& 5/7 \\
\PP( 1, 2, 2, 3, 3, 3, 4, 5 ) & 3\times\frac{1}{3},\frac{1}{5}(1,4,2)& 2/5& 7/15 \\
\PP( 1, 2, 2, 3, 3, 3, 4, 7 ) & 2\times\frac{1}{3}, \frac{1}{7}(3,4,2)& 8/21& 31/63 \\
\PP( 1, 2, 2, 3, 3, 4, 5, 5 ) & \frac{1}{3}, \frac{1}{5}(1,4,2), \frac{1}{5}(2,3,2)& 1/3& 22/45 \\
\PP( 1, 2, 2, 3, 3, 4, 5, 7 ) & \frac{1}{5}(2,3,2), \frac{1}{7}(3,4,2)& 11/35 & 18/35 \\
\hline
\PP( 1, 2, 2, 3, 3, 5, 8, 11 ) & \frac{1}{11}(3,8,2)& 3/11& 6/11 \\
\PP( 1, 2, 2, 3, 4, 5, 5, 7 ) & \frac{1}{5}(1,4,2), \frac{1}{7}(2,5,2)& 9/35& 18/35 \\
\PP( 1, 2, 3, 4, 4, 5, 5, 5 ) & 3\times\frac{1}{5}(1,4,2)& 1/5& 2/5 \\
\PP( 1, 2, 3, 4, 4, 5, 5, 9 ) & \frac{1}{5}(1,4,2), \frac{1}{9}(4,5,2)& 8/45& 58/135 \\
\PP( 1, 2, 3, 4, 4, 5, 9, 13 ) & \frac{1}{13}(4,9,2)& 2/13& 6/13 \\
\hline
\PP( 1, 2, 3, 4, 5, 5, 6, 7 ) & \frac{1}{3}, \frac{1}{5}(1,4,2), \frac{1}{7}(1,6,2)& 17/105& 127/315 \\
\PP( 1, 2, 3, 4, 5, 5, 6, 11 ) & \frac{1}{3}, \frac{1}{11}(5,6,2)& 5/33& 43/99 \\
\PP( 1, 2, 3, 4, 5, 6, 7, 7 ) & \frac{1}{7}(1,6,2), \frac{1}{7}(3,4,2)& 1/7& 3/7 \\
\PP( 1, 2, 3, 5, 6, 7, 8, 9 ) & 2\times\frac{1}{3}, \frac{1}{9}(1,10,2)& 1/9& 11/27 \\
\PP( 1, 2, 5, 7, 8, 9, 10, 11 ) & \frac{1}{5}(2,3,2), \frac{1}{11}(1,10,2)& 3/55& 19/55 \\
\hline
\PP( 2, 2, 3, 3, 4, 5, 5, 5 ) & \frac{1}{3},3\times\frac{1}{5}(2,3,2)& 2/15& 13/45 \\
\PP( 2, 2, 3, 3, 4, 5, 5, 7 ) & 2\times\frac{1}{3}, \frac{1}{5}(2,3,2), \frac{1}{7}(2,5,2) &13/105 & 92/315 \\
\PP( 2, 2, 3, 3, 4, 5, 7, 9 ) & 3\times\frac{1}{3}, \frac{1}{9}(2,7,2)& 1/9& 8/27 \\
\PP( 2, 2, 3, 5, 5, 7, 12, 17 ) & \frac{1}{17}(5,12,2)& 1/17& 5/17 \\
\PP( 2, 3, 3, 4, 5, 5, 6, 7 ) & 5\times\frac{1}{3},\frac{1}{5}(1,4,2)& 1/15& 11/45 \\
\hline
\PP( 2, 3, 3, 4, 5, 7, 10, 13 ) & 2\times\frac{1}{3}, \frac{1}{13}(3,10,2)& 2/39& 28/117 \\
\PP( 2, 3, 4, 5, 5, 6, 7, 7 ) & \frac{1}{3},\frac{1}{5}(1,4,2),\frac{1}{5}(2,3,2),\frac{1}{7}(3,4,2)& 1/21& 64/315 \\
\PP( 2, 3, 4, 5, 5, 6, 7, 9 ) & 2\times\frac{1}{3}, \frac{1}{5}(2,3,2), \frac{1}{9}(4,5,2)& 2/45& 28/135 \\
\PP( 2, 3, 5, 6, 7, 7, 8, 9 ) & 3\times\frac{1}{3}, \frac{1}{5}(2,3,2), \frac{1}{7}(1,6,2)& 1/35& 19/105 \\
\PP( 2, 5, 5, 6, 7, 8, 9, 11 ) & 2\times\frac{1}{5}(2,3,2), \frac{1}{11}(5,6,2)& 1/55& 8/55 \\
\hline
\end{array}
$
\caption{Fano 3-folds in codimension 4\label{tab!cod4}}
\end{table}

\vspace{5mm}
\noindent
Gavin Brown, IMSAS, University of Kent, CT2 7AF, UK.\\
Email: gdb@kent.ac.uk

\vspace{5mm}
\noindent
Kaori Suzuki, Tokyo Institute of Technology, 2-12-1 Ookayama, Meguro-ku,
152-8550, Japan.\\
Email:  k-suzuki@math.titech.ac.jp

\end{document}